%% file: main.tex
\documentclass[12pt]{article}
\usepackage[margin=1.25in]{geometry}
\usepackage{amsmath,amsthm,amssymb}
\usepackage{graphicx}
\usepackage{wasysym}
\usepackage{todonotes}
\usepackage[linesnumbered,ruled,vlined]{algorithm2e}
\usepackage{authblk}
\usepackage{subfiles}

\newtheorem{quest}{Question}{\bfseries}{\rmfamily}
\newtheorem{definition}{Definition}[section]
\newtheorem{theorem}{Theorem}[section]
\newtheorem{conjecture}{Conjecture}

\newtheorem{corollary}{Corollary}[theorem]
\newtheorem{lemma}[theorem]{Lemma}
\newcommand{\floor}[1]{\lfloor #1 \rfloor}
\newcommand{\ceil}[1]{\left\lceil #1 \right\rceil}
\newcommand{\gen}[1]{\ensuremath{\langle #1\rangle}}
\usepackage{xcolor}
\definecolor{light-gray}{gray}{0.95}
\newcommand{\code}[1]{\colorbox{light-gray}{\texttt{#1}}}
\usepackage{array}
\newcolumntype{L}{>{$}l<{$}}

\global\long\def\Aut{\operatorname{Aut}}

\title{Local Configurations in Union-Closed Families}
\author{Jonad Pulaj and Kenan Wood\thanks{Department of Mathematics and Computer Science, Davidson College, Davidson, NC 28036, \texttt{\{jopulaj, kewood\}@davidson.edu}}}
\date{}

\begin{document}

\maketitle

\subfile{sections/0-abstract}

\subfile{sections/1-intro}

\subfile{sections/3-FC-numbers}

\subfile{sections/2-symmetry}

\subfile{sections/4-poonen}

\subfile{sections/5-conclusion}

\bibliographystyle{plain}
\bibliography{bibliography}

\end{document}

%% file: sections/0-abstract.tex
\begin{abstract}
    The Frankl or Union-Closed Sets conjecture states that for any finite union-closed family of sets $\mathcal{F}$ containing some nonempty set, there is some element $i$ in the ground set $U(\mathcal F) := \bigcup_{S \in \mathcal{F}} S$ of $\mathcal{F}$ such that $i$ is in at least half of the sets in $\mathcal{F}$. In this work,
    we find new values and bounds for the least integer $FC(k, n)$ such that any union-closed family containing $FC(k, n)$ distinct $k$-sets of an $n$-set $X$ satisfies Frankl's conjecture with an element of $X$.
    Additionally, we answer an older question of Vaughan regarding symmetry in union-closed families and we give a proof of a recent question posed by Ellis, Ivan and Leader. 
    Finally, we introduce novel local configuration criteria through a generalization of Poonen's Theorem to prove the conjecture for many, previously unknown classes of families.
\end{abstract}

%% file: sections/1-intro.tex
\section{Introduction}\label{Sec:Intro}
Frankl's or the Union-Closed Sets conjecture is an open, well-known problem in extremal set theory. A finite family of finite sets $\mathcal{F}$ is \emph{union-closed} if for every $A, B \in \mathcal{F}$, it follows that $A \cup B \in \mathcal{F}$. Frankl's conjecture states that for any union-closed family $\mathcal{F}$ containing some nonempty set, there is some element $i$ in the ground set, or universe, of $\mathcal{F}$ defined as $U(\mathcal F) := \bigcup_{S \in \mathcal{F}} S$ such that $i$ is in at least half of the sets in $\mathcal{F}$.

Of the well-known techniques to tackle Frankl's conjecture, this work is concerned with the approach of local configurations~\cite{bruhn}, a method that aims to prove the conjecture for any union-closed family $\mathcal{F}$ satisfying some local conditions with respect to some fixed ground set $X \subseteq U(\mathcal{F})$. We believe recent developments~\cite{pulaj:2019,pulaj:2021}, including the work presented here, provide a new impetus into this line of research and its implications for Frankl's conjecture. In addition, some questions related to local configurations may be of independent interest since they are not implied by Frankl's conjecture.

The genesis of local configurations began with the well-known observations that any union-closed family containing a 1-set or a 2-set satisfies Frankl's conjecture with an element from the 1-set or 2-set (where a $k$-set is a set with $k$ elements). Poonen \cite{poonen} provided a complete characterization of families $\mathcal{A}$ such that every union-closed family containing $\mathcal{A}$ satisfies Frankl's conjecture with an element from $U(\mathcal{A})$. Such families $\mathcal{A}$ are called Frankl-Complete (FC). Families $\mathcal{A}$ that are not FC, are called Non-FC. As a consequence he showed that there is a union-closed family $\mathcal{F}$ containing a 3-set $A$ such that every element of $A$ is in strictly less than half the sets of $\mathcal{F}$; that is, $\{A\}$ is Non-FC. Using Poonen's Theorem and machine-assisted techniques, Morris \cite{morris} and Vaughan \cite{vaughan:2002} were able to characterize many FC-families on at most six elements. More recently, Pulaj \cite{pulaj:2019} exhibited the first efficient algorithm to completely characterize FC-families on at most 10-elements.

For a positive integer $k$, we define $[k] = \{1, \dots, k\}$ and let $\mathcal{P}([k])$ denote the power set of $[k]$. For $3 \leq k < n$, define $FC(k, n)$ to be the least integer $m$ such that any $\mathcal{A} \subseteq \mathcal{P}([n])$ containing $m$ distinct $k$-sets is FC. Morris \cite{morris} proved that $FC(3, n) \geq \floor{n/2} + 1$ for all $n \geq 4$ and conjectured that equality holds. Pulaj \cite{pulaj:2021} then proved Morris's conjecture showing $FC(3, n) = \floor{n/2} + 1$ for all $n \geq 4$. Morris also proved that $FC(4, 5) = 5$ and $7 \leq FC(4, 6) \leq 8$, while Marić, Vučković and Živković \cite{fc6_classification} provided a complete classification of all FC-families on six elements, showing $FC(4, 6) = 7$.

Our main contributions in this work are as follows. 
First, we algorithmically show:
\begin{itemize}
    \item $FC(4, 7) = 10, FC(4, 8) = 12, FC(5, 7) = 14, FC(6, 8) = 26$;
    \item $FC(4, 9) \geq 14, FC(4, 10) \geq 16, FC(5, 8) \geq 21, FC(5, 9) \geq 31, FC(5, 10) \geq 44, FC(6, 9) \geq 42, FC(6, 10) \geq 71, FC(7, 10) \geq 85$.
\end{itemize}
We also prove the following new upper bounds for general $n$, which follow from the above exact characterizations of small $FC(k, n)$ computations.
\begin{itemize}
    \item $FC(4, n) \leq 1 + \ceil{\frac{11}{1680} \cdot n(n-1)(n-2)(n-3)}$ for $n > 8$;
    \item $FC(5, n) \leq 1 + \ceil{\frac{13}{2520} \cdot n(n-1)(n-2)(n-3)(n-4)}$ for $n > 7$;
    \item $FC(6, n) \leq 1 + \ceil{\frac{5}{4032} \cdot n(n-1)(n-2)(n-3)(n-4)(n-5)}$ for $n > 8$.
\end{itemize}
As a consequence, we obtain $FC(4, 9) \leq 21$, $FC(5, 8) \leq 36$, and $FC(6, 9) \leq 76$. In contrast to previous works~\cite{pulaj:2019, pulaj:2021, eifler} where exact integer programming is used for verification of computational results, in our current work we use a SMT (Satisfiability Modulo Theory) solver for verification as suggested in~\cite{fc6_classification}. Tools like SMTCoq~\cite{smtcoq} pave the way for further verification in interactive theorem provers.

Second, we answer a question of Vaughan \cite{vaughan:2002} in the positive that simplifies Poonen's characterization of FC-families according to the symmetry of a given family, in particular, its automorphism group. For families $\mathcal{A}$ and $\mathcal{B}$, define $\mathcal{A} \uplus \mathcal{B} := \{A \cup B: A \in \mathcal{A}, B \in \mathcal{B}\}$. For an element $i$, let $\mathcal{A}_i := \{A \in \mathcal{A}: i \in A\}$. We explicitly state Poonen's Theorem below.
\begin{theorem}[Poonen]
\label{poonen}
Let $\mathcal{A}$ be a union-closed family of sets with $\emptyset \in \mathcal{A}$ and $U(\mathcal{A}) = [n]$. Then the following are equivalent:
\begin{enumerate}
    \item $\mathcal{A}$ is an FC-family. That is, for all union-closed $\mathcal{F} \supseteq \mathcal{A}$, there is some $i \in U(\mathcal{A})$ such that $|\mathcal{F}_i| \geq |\mathcal{F}|/2$.
    \item There exists some $c \in \mathbb{R}_{\geq 0}^n$ satisfying $\sum_{i\in [n]} c_i = 1$ such that for any union-closed $\mathcal{B} \subseteq \mathcal{P}([n])$ with $\mathcal{A} \uplus \mathcal{B} = \mathcal{B}$, we have
    \[
    \sum_{i \in [n]} c_i|\mathcal{B}_i| \geq |\mathcal{B}|/2.
    \]
\end{enumerate}
\end{theorem}
The set of all $c$ allowed in (2) is a polyhedron denoted $P^\mathcal{A}$. In 2002, Vaughan \cite{vaughan:2002} asked whether or not $P^\mathcal{A} \neq \emptyset$ implies there is some $c \in P^\mathcal{A}$ such that $c_i = c_j$ whenever there is an automorphism of $\mathcal{A}$ mapping $i$ to $j$. We prove that this implication does, indeed, hold.

Third, we highlight the utility of FC-families by answering in the positive the following recently posed question by Ellis, Ivan and Leader \cite{torus_conjecture}. Let $n \geq 4$ and choose some $R \subset \mathbb{Z}_n$ with $|R| = 3$. Does the union-closed family generated by all translates of $R \times \{0\}$ or $\{0\} \times R$ by elements of $\mathbb{Z}_n \times \mathbb{Z}_n$ necessarily satisfy Frankl's conjecture?

Finally, we prove a useful generalization of Poonen's Theorem that constructs a new type of local configuration. Given a Non-FC family $\mathcal{A}$, our theorem gives a method of restricting the possible union-closed families $\mathcal{F} \supseteq \mathcal{A}$ such that $|\mathcal{F}_i| < |\mathcal{F}|/2$ for all $i \in U(\mathcal{A})$ by considering only the structure of $\mathcal{A}$. To our knowledge, this is the first result that allows us to prove that a large collection of union-closed families that contain a \textit{Non-FC} family satisfy Frankl's conjecture. We are able to obtain very strong results about union-closed families containing a small Non-FC family such as $\{\emptyset, \{1,2,3\}\}$ or $\{\emptyset, \{1,2,3,4\}\}$. We obtain similar results on families of 4-sets, 5-sets, 6-sets, and 7-sets using an adaptation of Pulaj's algorithm together with the algorithm used for the above computations on $FC(k, n)$.

The rest of this paper is organized as follows.
In Section \ref{Sec:FC-numbers}, we give many new values and bounds for $FC(k, n)$, along with an interesting structural conjecture.
Section \ref{Sec:Symmetry} settles an older question of Vaughan, and Section \ref{Sec:Leader} settles a more recent question of Ellis, Ivan and Leader. 
Section \ref{Sec:Poonen} proves a generalization of Poonen's Theorem and introduces a new kind of local configuration, which we use to prove Frankl's conjecture for many new previously unknown classes of families.
Finally, concluding remarks may be found in Section \ref{Sec:Conclusion}.

%% file: sections/3-FC-numbers.tex
\section{FC-values and FC-bounds}\label{Sec:FC-numbers}
In this section, we give many new values and bounds implying Frankl-Completeness and conjecture a striking structural pattern regarding maximal Non-FC families.
\begin{definition}
Two families of sets $\mathcal{A}$ and $\mathcal{B}$ are \emph{isomorphic}, written $\mathcal{A} \cong \mathcal{B}$, provided there is some bijection $\phi: U(\mathcal{A}) \to U(\mathcal{B})$ such that $\mathcal{B} = \{\phi(S): S \in \mathcal{A}\}$. The map $\phi$ is called an \emph{isomorphism}.
\end{definition}

We now introduce our main tool for determining exact values of $FC(k, n)$, Algorithm 1. The method \code{getNonIsomorphicFamilies($n, k, m$)} returns a set of representatives from each isomorphism class of families $\mathcal{A}$ of $m$ distinct $k$-sets with $U(\mathcal{A}) = [n]$. The method \code{isFC($\mathcal{F}$)} uses Pulaj's algorithm \cite{pulaj:2019} to return true if $\mathcal{F}$ is FC, and false otherwise.

\SetKwInput{KwInput}{Input}
\SetKwInput{KwOutput}{Output}

\begin{algorithm}[!ht]
\DontPrintSemicolon
\caption{getNFC($n, k, m$)}
  
\KwInput{Positive integers $n, k, m$ where $n \geq k \geq 3$}
\KwOutput{A set of all pair-wise nonisomorphic Non-FC families of $m$ distinct $k$-sets with universe $[n]$}

\If{$km < n$ \emph{\textbf{or}} $m > \binom{n}{k}$}{
    \Return $\emptyset$
}

\;

$\mathbb{NFC} \gets \emptyset$

$\mathbb{FC} \gets \emptyset$

\;

\If{$k(m-1) < n$}
    {
    \For{$\mathcal F \in \emph{getNonIsomorphicFamilies}(n, k, m)$}
        {
        \If{\emph{\textbf{not}} $\emph{isFC}(\mathcal F)$}
            {$\mathbb{NFC} \gets \mathbb{NFC} \cup \{\mathcal F\}$}
        }
    \Return $\mathbb{NFC}$
    }

\;

$J \gets \{i \in \mathbb Z \mid \max\{k, n-k\} \leq i \leq n\}$

\For{$\mathcal{F} \in \bigcup_{i \in J}$ \emph{getNFC}$(i, k, m-1)$}
    {
    \For{$S \subseteq [n]$ such that $|S| = k$ and $U(\mathcal{F} \cup \{S\}) = [n]$}
        {
        \If{$\forall\mathcal{A} \in \mathbb{NFC} \cup \mathbb{FC}\colon \mathcal{F} \cup \{S\} \not\cong \mathcal{A}$}
            {
            \If{$\mathcal{F} \cup \{S\}$ contains a proper FC-family}{\textbf{continue}}
            
            \If{\emph{isFC}$(\mathcal{F} \cup \{S\})$}
                {$\mathbb{FC} \gets \mathbb{FC} \cup \{\mathcal{F} \cup \{S\}\}$}
            \Else
                {$\mathbb{NFC} \gets \mathbb{NFC} \cup \{\mathcal{F} \cup \{S\}\}$}
            }
        }
    }
\Return $\mathbb{NFC}$
\end{algorithm}

Algorithm 1 is a recursive algorithm designed to determine all isomorphism classes of Non-FC families of $m$ distinct $k$-sets with universe $[n]$, while disregarding families containing a proper FC-family. The isomorphism checks in line 16 are performed by computing a canonical form\footnote{We use SageMath's \code{canonical\_label()} method within the \code{IncidenceStructure} class.} of the family $\mathcal{F}$ such that any family isomorphic to $\mathcal{F}$ has an identical canonical form, checking if that canonical form has been computed before, and if not, storing its canonical form. The proper FC-containment check in line 17 is computed in a similar fashion by computing the canonical form of subfamilies of $\mathcal{F}$ with one fewer member-set. In our implementation, we manually start at the bottom of the call stack to avoid recomputation.

Additionally, for the purpose of ensuring the correctness of each \code{isFC()} computation, we use the SMT solver Z3 \cite{z3} within the SMT python library, pySMT \cite{pysmt}.
For verifying Non-FC families, we check the infeasibility of the terminating set of constraints produced by the \code{isFC()} algorithm. For FC families, (using Pulaj's notation) we check the infeasibility of the linear integer system defining $X(\mathcal{A}, c)$, where $c$ is the vector in $\mathbb Z^n$ found by the algorithm that is proposed to satisfy $X(\mathcal{A}, c) = \emptyset$.

\begin{lemma}
Algorithm 1 correctly finds a desired collection of Non-FC families.
\end{lemma}
\begin{proof}
For termination, notice that in each recursive call, we must have $n \geq k > 0$. Also, the $m$ argument is decreased by 1 at every call, so if Algorithm 1 did not terminate, $km \geq n$ at every iteration; however, $m$ must be zero at some point assuming no termination. This is a contradiction because $n > 0$. Therefore Algorithm 1 terminates.

For correctness, observe that the theorem is true if either $km < n$ or $k(m-1) < n$ or $m > \binom{n}{k}$. We first prove the following claim.\\

\noindent \underline{Claim:} Let $J = \{i \in \mathbb Z \mid \max\{k, n-k\} \leq i \leq n\}$. Assume getNFC($i, k, m-1$) is correct for all $i \in J$. Then getNFC($n, k, m$) is correct.\\

\noindent \textit{Proof of claim.} We may assume $k(m-1) \geq n$ and $m \leq \binom{n}{k}$. Consider the execution of getNFC($n, k, m$). Observe that anytime a family is added to $\mathbb{NFC}$, we always first verify that it is Non-FC. Hence every family in $\mathbb{NFC}$ is Non-FC. Suppose $\mathcal{F}$ is a Non-FC family with universe $[n]$ containing $m$ distinct $k$-sets. Let $S \in \mathcal{F}$; let $\mathcal{F}' = \mathcal{F}-\{S\}$ with $i := |U(\mathcal{F})|$. Since $S \in \mathcal{F}$, we know $|S| = k$, so that $n-k \leq i$ and $k \leq i \leq n$ (because $m \geq 2$). Hence $i \in J$, which implies that getNFC($n, k, m$) iterates through all families in getNFC($i, k, m-1$). By assumption, one of these families, say $\mathcal{G}'$, is isomorphic to $\mathcal{F}'$. Since there is an isomorphism $\phi: U(\mathcal{F}') \to U(\mathcal{G}')$, the family $\mathcal{G} := \mathcal{G}' \cup \{\phi(S \cap U(\mathcal{F}')) \cup (S - U(\mathcal{F}'))\}$ is isomorphic to $\mathcal{F}$. Also $\mathcal{G}$ is added to $\mathbb{NFC}$ since $\mathcal{F} \cong \mathcal{G}$ is Non-FC, as desired. Thus getNFC($n, k, m$) is correct.\\

We proceed by induction on $n$. Note that $n \geq k$, so the base case is $n = k$. If $m = 1$, the the result follows by inspecting lines 7-11 in Algorithm 1. If $m \geq 2$, then $m > \binom{n}{k} = 1$, so getNFC($n, k, m$) correctly returns.

For the induction step on $n$, suppose $n \geq k+1$ and getNFC($n', k, m'$) is correct for all $k \leq n' < n$ and $m'$. Then getNFC($i, k, m-1$) correctly returns for all $i \in J-\{n\}$. To show getNFC($n, k, m$) correctly returns, we use induction on $m$. If $m = 1$, then certainly getNFC($n, k, m$) is correct. Suppose $m \geq 2$ and getNFC($n, k, m-1$) correctly returns. This shows that getNFC($i, k, m-1$) is correct for all $i \in J$. Hence the theorem follows from the above claim.
\end{proof}
Using Algorithm 1, which can be easily extended to determine the exact value of $FC(k, n)$ for small values of $k$ and $n$, we have determined the following.\footnote{All code used in this paper can be accessed here: https://github.com/KenanWood/Local-Configurations-in-Union-Closed-Families}
\begin{theorem}\label{FC(4,7)}
$FC(4, 7) = 10$.
\end{theorem}

\begin{theorem}\label{FC(4,8)}
$FC(4, 8) = 12$.
\end{theorem}

\begin{theorem}\label{FC(5,7)}
$FC(5, 7) = 14$.
\end{theorem}

\begin{theorem}\label{FC(6,8)}
$FC(6, 8) = 26$.
\end{theorem}

The system used to verify all of our results (including those in Section \ref{Sec:Poonen}) has an Intel Xeon Processor E5-2620 v4 with 16 cores, each running at 2.1GHz; the system has 128GB of memory and two NUMA nodes. Theorems \ref{FC(4,7)}, \ref{FC(5,7)}, \ref{FC(6,8)} have been verified within at most a couple hours, but Theorem \ref{FC(4,8)} took us more than 26 days to verify.

Let $\binom{S}{k}$ be the set of all $k$-subsets of a set $S$. Define a strict total order $<$, called the lexicographic order, on the set $\binom{[n]}{k}$ by $A < B$ if $\min(A \Delta B)  \in A$. In this order, for fixed $n$ and $k$ and for any $S \in \binom{[n]}{k}$, define $[S] := \{A \in \binom{[n]}{k} \mid A \leq S\}$. Let $\{S_i^{n,k}\}_{i \geq 1} = \binom{[n]}{k}$, where $S_i^{n,k} < S_j^{n,k}$ for all $1 \leq i < j$. If $n$ and $k$ are clear, we simply write $S_i$ instead of $S_i^{n,k}$.

The following conjecture seems to be very promising based off of our experimental results.

\begin{conjecture}\label{Conj1}
For fixed $n > k \geq 3$, if $[S_m]$ is an FC-family for some positive integer $m$ and has universe size $n$, then $FC(k, n) \leq m$.
\end{conjecture}

This conjecture has been verified for all $n > k \geq 3$ such that $FC(k, n)$ is known (it is trivial for $k=3$ and any $n \geq 4$); there is always a maximum Non-FC family of the form $[S_m]$ for some $m$. If the conjecture is true, then we could easily find all exact values of $FC(k, n)$ for $n \leq 10$; all we would need to do in that case is to find an integer $m$ such that $[S_m]$ is FC and $[S_{m-1}]$ is Non-FC, giving us a result of $FC(k, n) = m$. 

It can be shown using the above method that the following FC lower bounds are also exact values, assuming Conjecture \ref{Conj1}. Most of these have been verified to be tight bounds within pySMT, except lower bounds of $FC(k, 10)$.

\begin{theorem}
$FC(4, 9) \geq 14$, $FC(4, 10) \geq 16$, $FC(5, 8) \geq 21$, $FC(5, 9) \geq 31$, $FC(5, 10) \geq 44$, $FC(6, 9) \geq 42$, $FC(6, 10) \geq 71$, $FC(7, 10) \geq 85$. The remaining values of $FC(k, n)$ for $5 \leq k < n \leq 10$ are undefined.
\end{theorem}
\begin{proof}
For each pair $(k, n) \in \{ (4, 9), (4, 10), (5, 8), (5, 9), (5, 10), (6, 9), (6, 10), (7, 10) \}$, the family $[S_{m-1}^{n,k}]$ as defined above is Non-FC, where $m$ is the proposed lower bound. For the remaining pairs $(k, n)$ when $5 \leq k < n \leq 10$, we can easily show that the family $\binom{[n]}{k}$ is Non-FC.
\end{proof}

Assuming Conjecture \ref{Conj1} is true, Table \ref{tab:FC} shows a complete classification of FC-values for $(k, n) \in \{3, \dots, 7\} \times \{4, \dots, 10\}$, where no entry at $(k, n)$ indicates that $FC(k, n)$ is undefined.
\begin{table}
    \centering
    \begin{tabular}{c|c c c c c c c}
        $k \backslash n$ & 4 & 5 & 6 & 7 & 8 & 9 & 10\\
        \hline
        3 & 3 & 3 & 4 & 4 & 5 & 5 & 6\\
        4 &   & 5 & 7 & 10 & 12 & 14 & 16\\
        5 &   &   &   & 14 & 21 & 31 & 44\\
        6 &   &   &   &    & 26 & 42 & 71\\
        7 &   &   &   &    &    &    & 85\\
    \end{tabular}
    \caption{FC-values}
    \label{tab:FC}
\end{table}

To find upper bounds of $FC(k, n)$, we generalize and tighten a result of Morris \cite{morris}. Morris showed that $FC(4, n) \leq \frac{7}{360}n^4$, though without an explicit proof. The following theorem improves and generalizes this bound, which yields improved explicit upper bounds on $FC(k, n)$ for $4 \leq k \leq 6$.

\begin{theorem}\label{general_upper_bounds}
If $m_0 = FC(k, n_0) \leq \binom{n_0}{k}$, then
\[FC(k, n) \leq 1 + \ceil{\frac{(m_0-1)}{n_0 \cdots (n_0-k+1)} \cdot n \cdots (n-k+1)} \leq \binom{n}{k}\]
for all $n > n_0$.
\end{theorem}
\begin{proof}
Suppose $m_0 = FC(k, n_0) \leq \binom{n_0}{k}$. Let $n > n_0$ and $m := 1 + \ceil{(m_0-1) \cdot \frac{(n_0-k)!}{n_0!} \cdot \frac{n!}{(n-k)!}}$, noting that $m = 1 + \ceil{\frac{(m_0-1)}{n_0 \cdots (n_0-k+1)} \cdot n \cdots (n-k+1)}$. Since $n > n_0$, we know $\frac{(n_0-k)!}{n_0!} \cdot \frac{n!}{(n-k)!} > 1$. This implies
\begin{align*}
    m &\leq \ceil{1+ \left(\binom{n_0}{k}-1\right) \cdot \frac{(n_0-k)!}{n_0!} \cdot \frac{n!}{(n-k)!}}\\
    &= \ceil{1 + \binom{n_0}{k} \cdot \frac{(n_0-k)!}{n_0!} \cdot \frac{n!}{(n-k)!} - \frac{(n_0-k)!}{n_0!} \cdot \frac{n!}{(n-k)!}}\\
    &\leq \ceil{\frac{n_0!}{k!(n_0-k)!} \cdot \frac{(n_0-k)!}{n_0!} \cdot \frac{n!}{(n-k)!}}\\
    &= \binom{n}{k}.
\end{align*}

Next, let $\mathcal{A}$ be a family of $m$ distinct $k$-sets with a universe of size at most $n$. Define $\mathcal{A}^0 := \mathcal{A}$; for $i \geq 0$, recursively define $\mathcal{A}^{i+1} := \mathcal{A}^i$ if $|U(\mathcal{A}^i)| < n-i$, and otherwise, $\mathcal{A}^{i+1} := \mathcal{A}^i - \mathcal{A}^i_x$, where $x \in U(\mathcal{A}^i)$ minimizes $|\mathcal{A}^i_x|$. It follows that $|U(\mathcal{A}^i)| \leq n-i$ for all $i \geq 0$ by induction. Since $FC(k, n_0) = m_0$, it suffices to prove $|\mathcal{A}^{n-n_0}| \geq m_0$.

To this end, for any $i \geq 0$, the pigeonhole principle shows that there is some $x \in U(\mathcal{A}^i)$ such that $|\mathcal{A}^i_x| \leq \frac{k|\mathcal{A}^i|}{|U(\mathcal{A}^i)|}$. If $|U(\mathcal{A}^i)| = n-i$, then
\begin{align*}
    |\mathcal{A}^{i+1}| &\geq |\mathcal{A}^i| - \frac{k|\mathcal{A}^i|}{|U(\mathcal{A}^i)|}\\
    &= |\mathcal{A}^i|\left(1-\frac{k}{n-i}\right)\\
    &= |\mathcal{A}^i|\left( \frac{n-i-k}{n-i}\right).
\end{align*}
Otherwise, $|\mathcal{A}^{i+1}| = |\mathcal{A}^{i}|$ by construction, so that the inequality still holds since $\frac{n-i-k}{n-i} < 1$. In writing
\[
|\mathcal{A}^{n-n_0}| = |\mathcal{A}^0| \cdot \prod_{i=0}^{n-n_0-1}\frac{|\mathcal{A}^{i+1}|}{|\mathcal{A}^i|},
\]
we obtain
\begin{align*}
    |\mathcal{A}^{n-n_0}| &\geq m \cdot \prod_{i=0}^{n-n_0-1} \left(\frac{n-i-k}{n-i}\right)\\
    &\geq \left(1 + (m_0-1) \cdot \frac{(n_0-k)!}{n_0!} \cdot \frac{n!}{(n-k)!}\right) \cdot \left(\frac{(n-k)!}{n!}\cdot \frac{n_0!}{(n_0-k)!}\right)\\
    &> m_0 - 1,
\end{align*}
so that $|\mathcal{A}^{n-n_0}| \geq m_0$. Thus $\mathcal{A} \supseteq \mathcal{A}^{n-n_0}$ is an FC-family and the result follows.
\end{proof}

\begin{corollary}\label{Cor:FC-general-upper-bounds}
The following bounds hold:
\begin{itemize}
    \item $FC(4, n) \leq 1 + \ceil{\frac{11}{1680} \cdot n(n-1)(n-2)(n-3)}$ for $n > 8$;
    \item $FC(5, n) \leq 1 + \ceil{\frac{13}{2520} \cdot n(n-1)(n-2)(n-3)(n-4)}$ for $n > 7$;
    \item $FC(6, n) \leq 1 + \ceil{\frac{5}{4032} \cdot n(n-1)(n-2)(n-3)(n-4)(n-5)}$ for $n > 8$.
\end{itemize}
\end{corollary}
\begin{proof}
This is an immediate consequence of Theorem \ref{general_upper_bounds} along with Theorems \ref{FC(4,8)}, \ref{FC(5,7)}, and \ref{FC(6,8)}.
\end{proof}
\begin{corollary}\label{Cor:FC-upper-bounds}
$FC(4, 9) \leq 21$, $FC(5, 8) \leq 36$, and $FC(6, 9) \leq 76$.
\end{corollary}
\begin{proof}
This is an immediate consequence of Corollary \ref{Cor:FC-general-upper-bounds}.
\end{proof}

%% file: sections/2-symmetry.tex
\section{Symmetry in FC-families}\label{Sec:Symmetry}
In this section, we answer two previously unsolved questions regarding symmetry in union-closed families with respect to local configurations.

Given a union-closed family $\mathcal{A}$ containing $\emptyset$ with $U(\mathcal{A}) = [n]$, let $\mathbb{B}(\mathcal{A})$ be the set of all union-closed $\mathcal{B} \subseteq \mathcal{P}([n])$ such that $\mathcal{A} \uplus \mathcal{B} = \mathcal{B}$. Recall that $P^\mathcal{A} = \{c \in \mathbb{R}_{\geq 0}^n: \sum_{i \in [n]} c_i = 1 \wedge \forall \mathcal{B} \in \mathbb{B}(\mathcal{A}), \sum_{i \in [n]} c_i|\mathcal{B}_i| \geq |\mathcal{B}|/2\}$. Then By Poonen's Theorem \ref{poonen}, $\mathcal{A}$ is FC if and only if $P^\mathcal{A} \neq \emptyset $. As outlined in our introduction, the following is a generalization of Vaughan's \cite{vaughan:2002} question.
    
\begin{quest}
Given a union-closed family $\mathcal{A}$ containing $\emptyset$, if $P^\mathcal{A}$ is nonempty, then is there always some $c \in P^\mathcal{A}$ such that $c_i = c_j$ whenever there is an automorphism of $\mathcal{A}$ that maps $i$ to $j$?
\end{quest}

We prove that the answer is yes in the following theorem. First, let $\Aut(\mathcal{A})$ denote the set of all automorphisms of $\mathcal{A}$ and note that $\Aut(\mathcal{A})$ is a group under function composition.
\begin{theorem}\label{Thm:vaughan-conjecture}
Let $\mathcal{A}$ be a union-closed family containing $\emptyset$. If $P^\mathcal{A}$ is nonempty, then there is some $c \in P^\mathcal{A}$ such that $c_i = c_j$ whenever there is an automorphism of $\mathcal{A}$ that maps $i$ to $j$.
\end{theorem}
\begin{proof}
Without loss of generality, assume $U(\mathcal{A}) = [n]$. Suppose $x \in P^\mathcal{A}$. For any $\phi \in \Aut(\mathcal{A})$, we first show that $(x_{\phi(i)})_{i\in [n]} \in P^\mathcal{A}$; it suffices to show that for any $\mathcal{B} \in \mathbb{B}(\mathcal{A})$, we have $\sum_{i \in [n]} x_{\phi(i)} |\mathcal{B}_i| \geq |\mathcal{B}|/2$. Consider the image $\phi(\mathcal{B}) = \{\phi(S): S \in \mathcal{B}\}$. If $A' \in \mathcal{A}$ and $B' \in \phi(\mathcal{B})$, then there are $A \in \mathcal{A}$ and $B \in \mathcal{B}$ such that $A' = \phi(A)$ and $B' = \phi(B)$; the first holds since $\phi^{-1} \in \Aut(\mathcal{A})$, so that $A = \phi^{-1}(A') \in \mathcal{A}$, and the second is by construction of $\phi(\mathcal{B})$; this shows $A' \cup B' = \phi(A \cup B) \in \phi(\mathcal{B})$, which implies $\phi(\mathcal{B}) \in \mathbb{B}(\mathcal{A})$. Since $x \in P^\mathcal{A}$, then $\sum_{i \in [n]} x_i |\phi(\mathcal{B})_i| \geq |\phi(\mathcal{B})|/2$. Since $\phi$ is a bijection, $\sum_{i \in [n]} x_{\phi(i)} |\phi(\mathcal{B})_{\phi(i)}| \geq |\phi(\mathcal{B})|/2$, which shows $\sum_{i \in [n]} x_{\phi(i)} |\mathcal{B}_i| \geq |\mathcal{B}|/2$. Thus $(x_{\phi(i)})_{i \in [n]} \in P^\mathcal{A}$ for any $\phi \in \Aut(\mathcal{A})$.

Consider the convex combination of elements of $P^\mathcal{A}$,
\begin{align*}
    c &= \frac{1}{|\Aut(\mathcal{A})|} \cdot \sum_{\phi \in \Aut(\mathcal{A})} (x_{\phi(i)})_{i\in [n]}\\
    &= \frac{1}{|\Aut(\mathcal{A})|} \cdot \left(\sum_{\phi \in \Aut(\mathcal{A})} x_{\phi(i)}\right)_{i\in [n]}.
\end{align*}
For each $i \in [n]$, let 
\[s_i = \sum_{\phi \in \Aut(\mathcal{A})} x_{\phi(i)}.\]
Then we obtain
\[
c = \frac{1}{|\Aut(\mathcal{A})|} \cdot (s_i)_{i\in [n]}.
\]
For any $\phi_0 \in \Aut(\mathcal{A})$, every automorphism in $\Aut(\mathcal{A})$ can be written as a unique left composition with $\phi_0$; that is, $\Aut(\mathcal{A}) = \{\phi \circ \phi_0 : \phi \in \Aut(\mathcal{A})\}$.
Then, for every $i, j \in [n]$ such that there is some automorphism $\phi_0 \in \Aut(\mathcal{A})$ mapping $i$ to $j$, we know 
\[
s_i = \sum_{\phi \in \Aut(\mathcal{A})} x_{[\phi \circ \phi_0](i)} = \sum_{\phi \in \Aut(\mathcal{A})} x_{\phi(j)} = s_j,
\]
so that $c_i = c_j$ as well. Since $c$ is a convex combination of points in a polyhedron, $c \in P^\mathcal{A}$.
\end{proof}

An important case of Theorem \ref{Thm:vaughan-conjecture} is for transitive families, where a family $\mathcal{A}$ is said to be \emph{transitive} if for any two elements $i,j \in U(\mathcal{A})$, there is an automorphism of $\mathcal{A}$ mapping $i$ to $j$. In particular, let $\mathcal{A}$ be a transitive union-closed family with universe $[n]$. Then Theorem \ref{Thm:vaughan-conjecture} implies that $\mathcal{A}$ is FC if and only if $(1/n)_{i \in [n]} \in P^\mathcal{A}$.

We also remark a computational consequence of this theorem as follows. Note that the dimension of $P^{\mathcal{A}}$ is at most $n-1$, and equality holds most of the time. However, through Theorem \ref{Thm:vaughan-conjecture}, $P^\mathcal{A}$ is nonempty if and only if the polyhedron obtained by adding the constraints from Theorem \ref{Thm:vaughan-conjecture} is also nonempty. But this constructed polyhedron has dimension at most $|[n]/\Aut(\mathcal{A})| - 1$, where $X/G$ denotes the set of orbits of elements in $X$ under the group action of $G$. It follows that for families of sets $\mathcal{A}$ where the number of distinct automorphism orbits is small (highly symmetric families), determining if some $c \in P^\mathcal{A}$ exists becomes much more computationally efficient.

\section{A Result on Transitive Families of 3-sets}\label{Sec:Leader}
In an Abelian group $(G, +)$, if $R \subseteq G$ and $g \in G$, we define the \emph{translation of $R$ by $g$ in $G$} as the set
\[
g + R := \{g + r : r \in R\}.
\]
The set of all translations (by some element $g \in G$) of $R$ in $G$ is denoted $T(R)$. Given a family $\mathcal{A}$, the \emph{union-closure} of $\mathcal{A}$, or the family \emph{generated} by $\mathcal{A}$, is defined as the union-closed family $\gen{\mathcal{A}} := \{\bigcup_{S \in \mathcal{A}'} S: \mathcal{A}' \subseteq \mathcal{A}\}$; note that $\gen{\mathcal{A}}$ contains $\emptyset$.

The authors of \cite{torus_conjecture} pose the following open question related to small sets in union-closed families.
\begin{quest}
    Given some 3-set $R \subset \mathbb{Z}_n$, does the union-closed family generated by $\mathcal{A} = T(R \times \{0\}) \cup T(\{0\} \times R) \subseteq \mathbb Z_n^2$ necessarily satisfy Frankl's conjecture?
\end{quest}
The authors remark that this family is transitive; that is, for any $x, y \in \mathbb{Z}_n^2$, there is an automorphism $\phi \in \Aut(\mathcal{A})$ such that $\phi(x) = y$.

Let $\mathcal{A}$ be a family of sets. Let $d(x) = |\mathcal{A}_x|$ be the \emph{degree} of $x$ in $\mathcal{A}$. The family $\mathcal{A}$ is \emph{regular} if $d(x) = d(y)$ for all $x,y \in U(\mathcal{A})$, in which case the \emph{degree} of $\mathcal{A}$ is the common degree.

\begin{lemma}\label{regularity}
Let $\mathcal{A}$ be a regular family of 3-sets with degree $k \geq 2$ and universe of size $n \geq 4$. Then $\mathcal{A}$ is FC. 
\end{lemma}
\begin{proof}
Since $\sum_{i \in U(\mathcal{A})} d(i) = \sum_{A \in \mathcal{A}} |A|$, we obtain $kn = 3m$, so that $m = kn/3 \geq 2n/3$, where $m = |\mathcal{A}|$. This implies that $m \ge \ceil{2n/3} \ge \floor{n/2} + 1 = FC(3, n)$, showing that $\mathcal{A}$ is FC.
\end{proof}
\begin{theorem}\label{Thm:leader}
Let $R \subset \mathbb{Z}_n$ be a 3-set, where $n \geq 4$. Then the family $\mathcal{A} = T(R \times \{0\}) \cup T(\{0\} \times R) \subseteq \mathbb Z_n^2$ is FC, and thus, $\gen{\mathcal{A}}$ satisfies Frankl's conjecture.
\end{theorem}
\begin{proof}
It is clear that $\mathcal{A}$ is a regular family of 3-sets with degree at least two. Furthermore, $|U(\mathcal{A})| \geq 4$. The result follows from Lemma \ref{regularity}.
\end{proof}

%% file: sections/4-poonen.tex
\section{A Generalization of Poonen's Theorem and FC-families}\label{Sec:Poonen}
In this section, we give a generalization of Poonen's Theorem that allows us to prove that a large class of union-closed families containing a potentially \emph{Non-FC} family $\mathcal{A}$ satisfies Frankl's conjecture with an element from $U(\mathcal{A})$. To our knowledge, this is the first result that allows us to obtain any significant information about union-closed families containing a \emph{Non-FC} family.

For fixed $n$, given a family $\mathcal{F}$ with $[n] \subseteq U(\mathcal{F})$ and a set $T \subseteq U(\mathcal{F}) \setminus [n]$, let $\mathcal{F}^{T, n} = \{S \cap [n]: S \in \mathcal{F}, S \setminus [n] = T\}$; if $n$ is clear, we simply write $\mathcal{F}^T$. For utility, observe that if $\mathcal{F}$ is union-closed containing $\mathcal{A}$ with $U(\mathcal{A}) = [n]$ and $\emptyset \in \mathcal{A}$, then for any $T \subseteq U(\mathcal{F}) \setminus [n]$, the family $\mathcal{F}^T$ is union-closed and $\mathcal{A} \uplus \mathcal{F}^T = \mathcal{F}^T$. This shows that condition 1 of Theorem \ref{Thm:poonen-general} is not vacuously true. Note that when $\mathbb{B}$ consists of all union-closed families $\mathcal{B} \subseteq \mathcal{P}([n])$ such that $\mathcal{A} \uplus \mathcal{B} = \mathcal{B}$, the statement of Theorem \ref{Thm:poonen-general} reduces precisely to Poonen's Theorem.
\begin{theorem}\label{Thm:poonen-general}
    Let $\mathcal{A}$ be a union-closed family containing $\emptyset$ with $U(\mathcal{A}) = [n]$. Let $\mathbb{B}$ be a set of union-closed subfamilies of $\mathcal{P}([n])$ such that for all $\mathcal{B} \in \mathbb{B}$, it follows that $\mathcal{A} \uplus \mathcal{B} = \mathcal{B}$. Assume $\mathcal{A} \in \mathbb{B}$. Then the following are equivalent:
    \begin{enumerate}
        \item For any union-closed $\mathcal{F} \supseteq \mathcal{A}$ where for any $T \subseteq U(\mathcal{F})\setminus [n]$, it follows that $\mathcal{F}^T$ is empty or $\mathcal{P}([n])$ or a family in $\mathbb{B}$, there is some $i \in [n]$ such that $|\mathcal{F}_i| \geq |\mathcal{F}|/2$.
    
        \item There is some $c \in \mathbb{R}_{\geq 0}^n$ where $\sum_{i \in [n]}c_i = 1$ and $\sum_{i \in [n]} c_i|\mathcal{B}_i| \geq |\mathcal{B}|/2$ for all $\mathcal{B} \in \mathbb{B}$.
    \end{enumerate}
\end{theorem}
\begin{proof}
    The proof follows that of Poonen's Theorem exactly, except instead of using \emph{all} union-closed families $\mathcal{B}$ such that $\mathcal{A} \uplus \mathcal{B} = \mathcal{B}$, we only use the families $\mathcal{B} \in \mathbb{B}$.
\end{proof}
Below we present a special case of Theorem \ref{Thm:poonen-general} that is easier to work with than the previous theorem. In particular, it allows for a simple extension of Pulaj's algorithm to determine if a family of sets $\mathcal{A}$ satisfies condition 1 below. 
\begin{theorem}
    Let $\mathcal{A}$ be a union-closed family containing $\emptyset$ with $U(\mathcal{A}) = [n]$. Let $\mathcal{V} \subseteq \mathcal{P}([n])$ with $\mathcal{A} \subseteq \mathcal{V}$. The following are equivalent:
    \begin{enumerate}
        \item For any union-closed family $\mathcal{F} \supseteq \mathcal{A}$ where for each $T \subseteq U(\mathcal{F})$, it follows that $\mathcal{F}^T$ is equal to $\mathcal{P}([n])$ or a subfamily of $\mathcal{V}$, there is some $i \in [n]$ such that $|\mathcal{F}_i| \geq |\mathcal{F}|/2$.

        \item There is some $c \in \mathbb{R}_{\geq 0}^n$ where $\sum_{i \in [n]}c_i = 1$ and $\sum_{i \in [n]} c_i|\mathcal{B}_i| \geq |\mathcal{B}|/2$ for all union-closed $\mathcal{B} \subseteq \mathcal{V}$ such that $\mathcal{A} \uplus \mathcal{B} = \mathcal{B}$.
    \end{enumerate}
\end{theorem}
\begin{proof}
    This is a consequence of Theorem \ref{Thm:poonen-general}, by choosing $\mathbb{B}$ to be the set of all union-closed $\mathcal{B} \subseteq \mathcal{V}$ such that $\mathcal{A} \uplus \mathcal{B} = \mathcal{B}$, noting that $\emptyset$ and $\mathcal{A}$ are families in $\mathbb{B}$.
\end{proof}
For brevity, any family $\mathcal{A}$ with $U(\mathcal{A}) = [n]$ together with a family $\mathcal{V} \subseteq \mathcal{P}([n])$ with $\mathcal{A} \subseteq \mathcal{V}$ is said to be \emph{$\mathcal{V}$-FC} if $\gen{\mathcal{A}}$ and $\mathcal{V}$ satisfy Theorem 5.2. That is, $\mathcal{A}$ is $\mathcal{V}$-FC if for any union-closed family $\mathcal{F} \supseteq \mathcal{A}$ where for each $T \subseteq U(\mathcal{F})$, it follows that $\mathcal{F}^T$ is equal to $\mathcal{P}([n])$ or a subfamily of $\mathcal{V}$, there is some $i \in [n]$ such that $|\mathcal{F}_i| \geq |\mathcal{F}|/2$.

Since most interesting cases are when $\mathcal{V}$ is union-closed, our implementation assumes that $\mathcal{V}$ is union-closed. This has the advantage of simply restricting the variables in the integer program $IP(\mathcal{A}, c)$ to only the ones indexed by sets in $\mathcal{V}$ instead of all sets in $\mathcal{P}([n])$, using Pulaj's \cite{pulaj:2019} notation $IP(\mathcal{A}, c)$. We cannot make this simple restriction if $\mathcal{V}$ is not union-closed because the union-closure inequalities in $IP(\mathcal{A}, c)$ require the variable indices to be closed under union. Alternatively, we could simply add the constraints $x_S = 0$ for all $S \in \mathcal{P}([n]) \setminus \mathcal{V}$. In either case, it is clear that making either of the above restrictions on the integer program $IP(\mathcal{A}, c)$ yields a correct algorithm for determining if $\mathcal{A}$ is $\mathcal{V}$-FC. We implement this algorithm in Gurobi \cite{gurobi}, and verify the results using pySMT in an analogous way to the $FC(k, n)$ results in Section \ref{Sec:FC-numbers}.

A natural candidate for a family $\mathcal{V}$ that will obtain strong results is $\mathcal{V} = \{S \in \mathcal{P}([n]): |S| \neq 1\}$. Through experimentation with Pulaj's \code{isFC()} algorithm, which iteratively finds the most restrictive inequalities in Poonen's Theorem, we find the following. For most families $\mathcal{A}$ in which we are able to determine 
\code{isFC($\mathcal{A}$)} on our system, the most restrictive inequalities in Poonen's Theorem are the ones induced by families of the form $\mathcal{B} = \mathcal{A} \uplus \mathcal{P}([n]\setminus \{i\})$ for $i \in [n]$. In fact, Morris \cite{morris} conjectured that these are the only inequalities needed. However, Pulaj \cite{pulaj:2019} disproved his conjecture. Still, it is reasonable to expect that removing these inequalities (by restricting $\mathcal{V}$ as above) will yield many new strong results about Non-FC families. In the following corollaries, we show that setting $\mathcal{V} = \{S \in \mathcal{P}([n]): |S| \neq 1\}$ does indeed produce striking new information about small Non-FC families.

\begin{corollary}
    Let $\mathcal{A} = \gen{\{\{1,2,3\}\}}$ and $\mathcal{V} = \mathcal{P}([3])\setminus\{\{1\}\}$. Then $\mathcal{A}$ is $\mathcal{V}$-FC.
\end{corollary}
\begin{corollary}
    Let $\mathcal{A} = \gen{\{\{1,2,3,4\}\}}$ and $\mathcal{V} = \{S \in \mathcal{P}([4]): |S| \neq 1\}$. Then $\mathcal{A}$ is $\mathcal{V}$-FC.
\end{corollary}
\begin{corollary}\label{FC_V(5,6)}
    Let $\mathcal{A} = \gen{\{\{1,2,3,4,5\}, \{1,2,3,4,6\}, \{1,2,3,5,6\}\}}$ and $\mathcal{V} = \{S \in \mathcal{P}([6]): |S| \neq 1\}$. Then $\mathcal{A}$ is $\mathcal{V}$-FC. However, if $\mathcal{A} = \gen{\{\{1,2,3,4,5\}, \{1,2,3,4,6\}\}}$ with the same $\mathcal{V}$, then $\mathcal{A}$ is not $\mathcal{V}$-FC.
\end{corollary}
Let $FC_\mathcal{V}(k, n)$ be the minimal $m$ such that any family $\mathcal{A} \subseteq \mathcal{P}([n])$ containing at least $m$ distinct $k$-sets and has $U(\mathcal{A}) = [n]$ is $\mathcal{V}$-FC. When $\mathcal{V} = \{S \in \mathcal{P}([6]): |S| \neq 1\}$, Corollary \ref{FC_V(5,6)} implies $FC_\mathcal{V}(5, 6) = 3$, which is a significant improvement on FC-families since $FC(5,6)$ is not even defined. To determine $FC_\mathcal{V}(k, n)$, we may simply adapt Algorithm 1, or for our purposes since the cases are sufficiently small, we simply run a brute-force check of all desired non-isomorphic families. To this end, we obtain the following.
\begin{corollary}
    Let $\mathcal{V} = \{S \in \mathcal{P}([7]): |S| \ne 1\}$. Then $FC_\mathcal{V}(5, 7) = 5$ and $FC_\mathcal{V}(6, 7) = 7$.
\end{corollary}

%% file: sections/5-conclusion.tex
\section{Conclusion}\label{Sec:Conclusion}
In this work, we answer two previously unsolved questions. One is an older question of Vaughan that shows the dimension of Poonen's polyhedron $P^\mathcal{A}$ can be reduced from $|U(\mathcal{A})|$ to the number of orbits of $\Aut(\mathcal{A})$ through a projection, which shrinks the search space and reduces computational work. We also answer a question of Ellis, Ivan and Leader related to union-closed families generated by 3-sets. Our solution highlights the continual importance of FC-families in that they provide simple solutions to difficult problems related to the Union-Closed Sets conjecture.
Furthermore, we find and verify many new values of $FC(k, n)$ for $k \geq 4$, of which only two were previously known. 
These computations lead to three new general upper bounds on $FC(4, n), FC(5, n), FC(6, n)$. Additionally, an insightful pattern emerges in the maximum Non-FC families through executions of Algorithm 1 (see Conjecture \ref{Conj1}), suggesting that Non-FC families have significantly more structure than previously thought.
Finally, we introduce a new class of local configurations, a kind of ``partial Frankl-Completeness,'' by generalizing Poonen's theorem. We use an adaptation of Pulaj's algorithm to obtain strong new results in this direction.

We believe several directions merit further attention, including Conjecture \ref{Conj1} and Theorem \ref{Thm:poonen-general}. In particular, what are the limits of Theorem \ref{Thm:poonen-general} and $\mathcal{V}$-FC families? FC-families were very useful for proving that Frankl's conjecture holds for all union-closed families $\mathcal{F}$ with $|U(\mathcal{F})| \leq 12$, so these newfound restrictions may be sufficient to prove Frankl's conjecture for larger ground sets $U(\mathcal{F})$, such as when $|U(\mathcal{F})|$ equals 13 or 14. How far can the results about $\mathcal{V}$-FC families be taken in relation to FC-families? Affirmative answers to these questions would be a significant step towards a deeper understanding of local configurations.